\documentclass[12pt]{article}

\usepackage{amsmath,cd2}
\usepackage{amsfonts,amssymb}

\newtheorem{satz}{Theorem}[section]
\newtheorem{defi}[satz]{Definition}

\newtheorem{kor}[satz]{Corollary}

\newtheorem{prop}[satz]{Proposition}

\newcounter{Roma}
\newcounter{Ara}
\newcounter{let}

\makeindex

\begin{document}

\newcommand{\nc}{\newcommand}

\nc{\mapco}{\,\colon\, }
\nc{\ab}{^{ab}}

\nc{\catc}{{\C}}
\nc{\we}{\vee}

\nc{\hra}{\hookrightarrow}

\nc{\epi}{epimorphism}
\nc{\repi}{regular epimorphism}
\nc{\mono}{monomorphism}
\nc{\iso}{isomorphism}
\nc{\coker}[1]{\mbox{${\rm Coker}(#1)$}}
\nc{\Ker}[1]{\mbox{Ker$(#1)$}}
\nc{\defgl}{\stackrel{def}{=}}

\nc{\V}{\vspace{3mm}}
\nc{\VV}{\vspace{4mm}}
\nc{\lra}{\longrightarrow}
\nc{\lla}{\longleftarrow}
\nc{\mr}[1]{ \stackrel{#1}{\lra} }
\nc{\ml}[1]{ \stackrel{#1}{\lla} }
\nc{\hmr}[1]{\hspace{2mm}\stackrel{#1}{\lra}\hspace{2mm}}
\nc{\hml}[1]{\hspace{2mm} \stackrel{#1}{\lla}\hspace{2mm}}
\nc{\N}{\noindent}
\nc{\st}{^{\prime}}
\nc{\ot}{\otimes}
\nc{\hcong}{ \hspace{2mm}\cong\hspace{2mm}  }
\nc{\hfbox}{\hfill$\Box$}
\nc{\REF}[1]{(\ref{#1})}

\def\Z{\ifmmode{Z\hskip -4.8pt Z} \else{\hbox{$Z\hskip -4.8pt Z$}}\fi}

\def\Q{\ifmmode{Q\hskip-5.0pt\vrule height6.0pt depth
0pt\hskip6pt}\else{\hbox{$Q\hskip-5.0pt\vrule height6.0pt depth
0pt\hskip6pt$}}\fi}

\newcommand{\NN}{\mbox{$I\!\!N$}}

\nc{\Ph}{\phantom{}}

\nc{\BE}{\begin{equation}}
\nc{\EE}{\end{equation}}

\nc{\dst}{\displaystyle}
\nc{\sst}{\scriptscriptstyle}
\nc{\ssst}{\scriptscriptstyle}
\nc{\proof}{\N{\bf Proof\,:}\quad}
\nc{\proofof}[1]{\N{\bf Proof of {#1}\,:}\quad}
\nc{\proofofthm}[1]{\N{\bf Proof of theorem \ref{#1}\,:}\quad}

\nc{\htt}[1]{^{\otimes #1}}

\newcommand{\ssur}[2]{\mbox{$#1 \!\to\!\!\!\!\!\to\! #2$}}

\def\mapsel#1{\mbox{$\rule[-5mm]{0mm}{12mm} \searrow\rlap
{$\vcenter{\makebox[0mm][r]{$\scriptstyle#1$\hspace{12.2mm}}}$}$}}

\def\surltop#1{\makebox[1cm]{\mbox{$\stackrel{#1}{\makebox[0mm]{$\lla$}\hspace{0.7mm}
\makebox[0mm]{$\lla$}}$}}}

\nc{\Sur}[1]{\mbox{$\:\stackrel{#1}{\lra\!\!\!\!\!\to\,}\:$}}

\def\INJ{\mbox{\mathsurround=0pt
\makebox[0mm][r]{\parbox{0mm}{\rule[-0.65mm]{0mm}{0.2mm}$\scriptscriptstyle>$}}
\makebox[0.7cm][l]{\parbox{0.7cm} {$\lra$}}}}

\def\Inj#1{\mbox{$\:\stackrel{#1}{\INJ}\:$}}

\nc{\Injup}[1]{\mapup{#1}}
\nc{\injup}[1]{\mapup{#1}}
\nc{\injdown}[1]{\mapdown{#1}}


\def\mapup#1{\mbox{$\rule[-1mm]{0cm}{0.7cm}  
\makebox[0mm][r]{\raisebox{0.2mm}{$\scriptstyle\phantom{\cong}$}\hspace{0.6mm}}
\bigg\uparrow\rlap{$\vcenter{\hbox{$\scriptstyle#1$}}$}$}}

\def\isoup#1{\mbox{$\rule[-1mm]{0cm}{0.7cm}  
\makebox[0mm][r]{\raisebox{0.2mm}{$\scriptstyle\cong$}\hspace{0.6mm}}
\bigg\uparrow\rlap{$\vcenter{\hbox{$\scriptstyle#1$}}$}$}}

\def\mapdown#1{\mbox{$\rule[-1mm]{0cm}{0.7cm}  
\makebox[0mm][r]{\raisebox{0.2mm}{$\scriptstyle $}\hspace{0.6mm}}
\bigg\downarrow\rlap{$\vcenter{\hbox{$\scriptstyle#1$}}$}$}}

\def\isodown#1{\mbox{$\rule[-1mm]{0cm}{0.7cm}  
\makebox[0mm][r]{\raisebox{0.2mm}{$\scriptstyle\cong$}\hspace{0.6mm}}
\bigg\downarrow\rlap{$\vcenter{\hbox{$\scriptstyle#1$}}$}$}}

\def\isor#1{\mbox{$\smash{\mathop{\longrightarrow}\limits^{\cong}_{#1}}$}}

\def\isol#1{\mbox{$\smash{\mathop{\longleftarrow}\limits^{\cong}_{#1}}$}}

\def\surdown#1{\makebox[0mm]{$\mapdown{\raisebox{0.4mm}{$\scriptstyle#1$}}$}
\makebox[0mm]{\raisebox{-2.15mm}{$\downarrow$}}}

\def\surup#1{\makebox[0mm]{$\mapup{\raisebox{0.4mm}{$\scriptstyle#1$}}$}
\makebox[0mm]{\raisebox{0.7mm}{$\mapup{}$}}}

\newcommand{\surr}[2]{@\cdhgeneric>->->\twoheadrightarrow>#1>#2>}

\newcommand{\surl}[2]{@\cdhgeneric>\twoheadleftarrow>->->#1>#2>}

\newcommand{\brokenr}[2]{
@\cdhgeneric>\raise2.8pt\hbox to3.5pt{\hrulefill}\mkern9mu>
\raise2.8pt\hbox to3.5pt{\hrulefill}\hbox to5pt{}>
\mkern-7.5mu\dashrightarrow>#1 >
#2>}

\newcommand{\brokenup}[1]{
@\cdvgeneric>\hat\cdot
>\raisebox{3pt}{$\vdots$}>
\vbox{\kern3.5pt\hbox{$\cdot$}\kern-3.5pt}
> >\hspace{1mm}#1>}

\newcommand{\functlr}[2]{
\raisebox{0.4pt}{$\hss\begin{CD}
@>\vbox{\hbox to 0pt{$\hss\begin{CD}@<#1<<\end{CD}\hss$}\vskip-2pt}
>#2 >
\end{CD}\hss$}
}

\newcommand{\maprr}[2]{
\raisebox{-0.9pt}{$\hss\begin{CD}
@>\vbox{\hbox to 0pt{$\hss\begin{CD}@>#1>>\end{CD}\hss$}\vskip-3pt}
>#2 >
\end{CD}\hss$}
}

\newcommand{\mapdd}[2]{
@\cdvstandard>\downarrow\hspace{1.5pt}\kern-4pt\hspace{1.5pt}\Big\downarrow>#1>#2>}

\newcommand{\mapud}[2]{
@\cdvstandard>\uparrow\kern-3pt\Big\downarrow>#1>#2>}

\newcommand{\sepi}[3]{\,\mbox{$#1\,: \ssur{#2}{#3}$}\,}

\nc{\auf}{\twoheadrightarrow}
 
\nc{\ruled}{\rule[-4mm]{0mm}{0mm}}

\nc{\qu}{quadratic}
\nc{\GBoGB}{G/BG\st \otimes G/BG\st}
\nc{\lstar}{_{\raisebox{-1mm}{$*$}}}
\nc{\sm}{\:{ \wedge}\:}

\nc{\Rmod}{${\bf R}$-module}

\nc{\map}[3]{\mbox{$#1 \mapco #2 \to #3$}}

\nc{\rond}{{\,\sst \circ\,}}

\nc{\ruleu}{\rule{0mm}{7mm}}

\nc{\T}[1]{\tilde{#1}}

\nc{\Imm}[1]{\mbox{${\rm Im}(#1)$}}

\nc{\tw}{\end{document}}

\nc{\QJG}{\frac{\dst I(G) J}{\dst I^2(G) J}}
\nc{\QJH}{\frac{\dst I(H) J}{\dst I^2(H) J}}
\nc{\otz}{\ot}
\nc{\UL}[2]{{\rm U}_{#1}{\rm L}(#2)}

\nc{\IRN}[1]{I_{R,\cal G}^{#1}(G)}

\nc{\ULG}[1]{{\rm U}_{#1}{\rm L}^{\cal G}(G)}
\nc{\UGH}[1]{{\rm U}_{#1}^{\cal G}(G,H)}

\nc{\LG}[1]{{\rm L}^{\cal G}_{#1}(G)}

\nc{\calG}{{\cal G}}
\nc{\AB}{^{AB}}

\title{On the second cohomology of semidirect products}
\author{Manfred Hartl and S\'ebastien Leroy}
\maketitle

\begin{center}
\N LAMAV, 
ISTV2, \\ Universit\'e de Valenciennes et du Hainaut-Cambr\'esis, \\
Le Mont Houy,  59313 Valenciennes Cedex 9, France.\\ 
Email: Manfred.Hartl@univ-valenciennes.fr\\ Phone n$^o$ 0033/327511901, Fax
0033/327511900.

\end{center}
\vspace{8mm}

\begin{abstract}
Let $G$ be a group which is the semidirect product of a normal subgroup $N$ and a subgroup $T$, and let $M$ be a $G$-module with not necessarily trivial $G$-action. Then we embed the simultaneous restriction map 
$res=(res^G_N,res^G_T)^t \mapco  H^2(G,M) \to H^2(N,M)^T \times H^2(T,M)$ into a natural five term exact sequence consisting of one and two-dimensional cohomology groups of the factors $N$ and $T$. The elements of $H^2(G,M)$ are represented in terms of group extensions of $G$ by $M$ constructed from extensions of $N$ and $T$.
\end{abstract}\vspace{8mm}

\nc{\Aut}{{\rm Aut}}
\nc{\AutM}{{\rm Aut}^M}
\nc{\tphi}{\T{\varphi}}

\nc{\EN}{\mbox{$\underline{E}_N\mapco M \Inj{i_N} E_N \Sur{\pi_N} N$}}
\nc{\ET}{\mbox{$\underline{E}_T\mapco M \Inj{i_T} E_T \Sur{\pi_T} T$}}
\nc{\Der}{\mbox{Der$(N,M)$}}

\N{\large \bf Introduction.}\quad The low dimensional cohomology groups $H^n(G,M)$, $n\le 2$, of a group $G$ with coefficients in a $G$-module $M$  crucially occur in many fields,  in algebra as well as in geometry. In fact, they reflect the structure of $G$ (and of $M$ if the $G$-action on it is non trivial) in a subtle way which is far from being understood in general. If $G$ admits a proper normal subgroup $N$ it can be viewed as an extension
\BE\label{NGQ} 1 \to N \to G \to Q \to 1\:, \EE
and one wishes to express the cohomology of $G$ in terms of the cohomology of the simpler ``pieces" $N$ and $Q$. Formally, the Lyndon-Hochschild-Serre spectral sequence (referred to as LHSSS in the sequel) $H^p(Q,H^{n-p}(N,M)) \Rightarrow H^n(G,M)$ solves this problem, computing certain filtration quotients of $H^n(G,M)$  provided one can manage to compute the corresponding differentials; those concerning $H^2(G,M)$ were determined by Huebschmann \cite{HueSS}, in terms of automorphism groups of group extensions and of 2-fold crossed extensions, data which, however, are not easy to control in general. Also, knowing the filtration quotients of $H^n(G,M)$ does not amount to knowing its group structure completely unless $M$ is a vector space, and one often needs to represent the elements of the abstract group $H^n(G,M)$ by either explicit cocycles or  group extensions (for $n=2$). Another approach to the study of $H^2(G,M)$ consists in embedding it into exact sequences involving the cohomology groups of $N$ and $Q$; the so-called ``fundamental exact sequence" derived from the LHSSS being
\[0 \to H^1(Q,M^N) \mr{inf} H^1(G,M) \mr{res} H^1(N,M)^Q \mr{d_2} H^2(Q,M^N) \mr{inf} H^2(G,M)_1\]
 \BE\label{7TES}  \hspace{-14.5mm} \mr{tr} H^1(Q,H^1(N,M)) \mr{d_2} H^3(Q,M^N) \EE
where $H^2(G,M)_1 = \Ker{res^G_N\mapco H^2(G,M)\to H^2(N,M)}$. (We remark that in \cite{HaRo} we offer an elementary conceptual construction and proof of this exact sequence which, unlike the one in  \cite{Hue8TS} concerning the first five terms,  does not invoke automorphism groups). A different extension of the first five terms of \REF{7TES}, embedding the full group $H^2(G,M)$ instead of only $H^2(G,M)_1$, is given by Huebschmann in \cite{Hue8TS}, as follows:
\BE\label{8TES} H^2(Q,M^N) \mr{inf} H^2(G,M) \mr{} {\rm Xpext}(G,N;M) \mr{\Delta} H^3(Q,M^N) \mr{inf} H^3(G,M) \EE
where the group ${\rm Xpext}(G,N;M)$ consists of equivalent classes of crossed pairs introduced in that paper. We point out that according to both sequences \REF{7TES} and \REF{8TES},
the study of  $H^2(G,M)$  involves a {\em three}-dimensional cohomology group of $Q$. 

If the extension \REF{NGQ}   {\em splits}\/, i.e.\ if $G$ is the semidirect product of $N$ and some subgroup $T$ isomorphic with $Q$, the situation is slightly better, at least in special cases. For example, if $T$ is free, the LHSSS amounts to a short exact sequence 
\BE\label{TfreeES} 0 \to H^1(T,H^1(N,M)) \to H^2(G,M) \mr{res} H^2(N,M)^Q \to 0 \EE
If $T$ is arbitrary, but $G$ acts {\em trivially}\/ on $M$, $H^2(G,M)$ contains $H^2(T,M)$ as a canonical direct factor, and the complementary piece
$H^2(G,M)_2 =  {\rm Ker}(res^G_T\mapco$ $ H^2(G,M)\to H^2(T,M))$
 fits into an exact sequence
\[0 \to H^1(T,{\rm Hom}(N,M)) \to H^2(G,M)_2 \mr{res} H^2(N,M)^T \to H^2(T,{\rm Hom}(N,M))\] 
\BE\label{TaES} \hspace{86mm}\to H^3(G,M)_2\EE
due to Tahara \cite{Ta} who also provides a construction of the elements of $H^2(G,M)_2$ in terms of suitable cocycles. Moreover, these results determine $H^2(G,M)$ out of only 1- and 2-dimensional cohomology groups of $N$ and $T$, in contrast with the sequences  \REF{7TES} and \REF{8TES}. When $M$ is a {\em non trivial}\/ $G$-module, however, $H^2(T,M)$ is no longer a direct factor of $H^2(G,M)$ if $M\neq M^N$, and no general description of the latter group in terms of first and second cohomology groups of $N$ and $T$ seems to be known.
This is now provided in the present paper, by embedding the ``simultaneous restriction map" 
\[ res=(res^G_N,res^G_T)^t \mapco  H^2(G,M) \hmr{} H^2(N,M)^T \times H^2(T,M) \]
into an exact sequence which generalizes both sequences \REF{TfreeES} and \REF{TaES}, as follows.
\nc{\shmr}[1]{\hspace{1mm}\mr{#1}\hspace{1mm}}
\[H^1(T,M) \shmr{\partial^0_{N*}} H^1(T,\Der) \shmr{\tau} H^2(G,M) \shmr{res} H^2(N,M)^T \times H^2(T,M) \]
\BE\label{H2sequ}  \hspace{93.3mm} \shmr{\phi} H^2(T,\Der) \EE
Here $\Der$ denotes the group of derivations from $N$ to $M$, which can be easily determined using Fox differential calculus, by means of the Jacobian associated to a presentation of $N$, see \cite{Brown}. Thus the two terms left of $H^2(G,M)$ are easily accessible to computation. 
The maps in sequence \REF{H2sequ} are described in theorem \ref{theorem} below. Note that our sequence, unlike the preceding ones, invokes the group  
$\Der$ instead of its quotient $H^1(N,M)$; this may be considered as the price to pay for avoiding the appearance of a cohomology group of dimension three.

As did Tahara in his work, we also construct the elements of $H^2(G,M)$ out of those of the other groups, but not in form of cocycles but of group extensions of $G$ by $M$, the basic idea being to somehow lift the semidirect product decomposition of $G$ to any  group $E$ fitting into an extension $0\to M \to E \mr{\pi} G \to 1$. In fact, $E$ turns out to be an ``amalgamated semidirect product" $E_N \rtimes_M E_T$ where $E_N=\pi^{-1}N$ and $E_T=\pi^{-1}T$; so sequence \REF{H2sequ} arises from studying the appropriate actions of $E_T$ on $E_N$, by using automorphisms of group extensions as Huebschmann did in the cited papers, but in a different way.

We also point out that a description of $H^2(G,M)$ for $G=N\rtimes T$ is given in \cite{Ga-Ma}, in terms of generators and relations computed from compatible presentations of $G$, $N$ and $T$.

\section{The main result}

Throughout this paper, $G$ denotes  a group and $M$ a $G$-module, i.e.\ an abelian group endowed with a not necessarily trivial action  $\psi\mapco G \to \Aut(M)$. As usual, $M^G$ denotes the subgroup of elements of $M$ which are invariant under the action of $G$.
Moreover, $(C^*(G,M),\partial_G^*)$ denotes the standard complex of normalized cochains on $G$ with values in $M$, i.e.\ $C^n(G,M)$ is the abelian group
of all functions $\beta \mapco G^{\times n} \to M$ annihilating any tuple $(g_1,\ldots,g_n)$ where $g_i=1$ for some $i$, and the differential $\partial_G^n \mapco C^n(G,M) \to C^{n+1}(G,M)$ is given by the formula
\[ \partial_G^n(\beta) (g_1,\ldots,g_{n+1}) = g_1 \beta(g_2,\ldots,g_{n+1}) + \sum_{i=1}^n (-1)^i \beta(g_1,\ldots,g_ig_{i+1}, \ldots,g_{n+1}) \]\[ \hspace{50mm} {} + (-1)^{n+1} \beta(g_1,\ldots,g_n)\]
By definition, $H^n(G,M) = H^n(C^*(G,M),\partial_G^*)$. Denote the group of $n$-cocycles of  $(C^*(G,M),\partial_G^*)$ by   $Z^n(G,M) =\Ker{\partial_G^n}$; in particular, ${\rm Der}(G,M) = \Ker{\partial^1_G}$ is the set of derivations from $G$ to $M$, i.e.\ the set of all functions $d\mapco G\to M$ such that
$d(gg')=gd(g') + d(g)$ for $g,g'\in G$.
If $N$ and $T$ are  subgroups of $G$ such that $N$ is normal then the action of $T$ on $N$ by conjugation, ${}^tn=tnt^{-1}$, induces an action  of $T$ on
  the complex  $(C^*(N,M),\partial_N^*)$, given by $(t\beta)(n_1,\ldots,n_k) = t\beta({}^{t^{-1}}n_1,\ldots, {}^{t^{-1}}n_k)$ for $\beta\in C^k(N,M)$. Thus $T$  acts on $H^*(N,M)$. The following elementary construction turns out to be crucial in the sequel.   If $\Gamma$ is a group then a homomorphism of $\Gamma$-modules $f\mapco N\to N'$ gives rise to the composite homomorphism
  \[ \omega^n_{\Gamma}(f) \mapco H^n(\Gamma, {\rm coker}(f)) \mr{\omega_n} H^{n+1}(\Gamma, {\rm Im}(f)) \mr{\omega_{n+1}}   H^{n+1}(\Gamma, \Ker{f}) \]
  where $\omega_n$ and $\omega_{n+1}$ are the connecting homomorphisms associated with the obvious short exact sequences of $\Gamma$-modules. In particular, for $\Gamma = T$ and $f=\hat{\partial}^1_N \mapco C^1(N,M) \to Z^2(N,M)$ given by   ${\partial}^1_N$ we get the map
  \[ \omega^0_T(\hat{\partial}^1_N) \mapco H^2(N,M)^T \to H^2(T,\Der) \:.\]
  The following conceptual construction of this map will be provided  
  in the proof of Proposition \ref{Autext}. Let $z \in H^2(N,M)^T$ be represented by a group extension $e\mapco M\Inj{} E \Sur{} N $ of $N$ by $M$, see section 2. Then 
  $\omega^0_T(\hat{\partial}^1_N)(z)$ is represented by the restriction to $T$ of the class of the extension
   \[ 0 \mr{}  \Der \mr{} \Aut_G(e) \mr{} G \mr{} 1 \]
   constructed by Huebschmann in \cite{Hue8TS}. More explicitly, relation \REF{Autclassexpl} below provides the following description of this class in terms of cocycles:
   Let $z$ be represented by a 2-cocycle $\beta\in C^2(N,M)$.  
  Then for $t\in T$ there exists $\gamma(t) \in C^1(N,M)$ such that 
 \BE\label{defgammaT} t\beta - \beta = {\partial}^1_N(\gamma(t) ) \:.\EE
 We thus get a map $\gamma \in C^1( T,C^1(N,M))$. Its image ${\partial}^1_{T}(\gamma) \in Z^2(  T,C^1(N,M))$ actually takes values in ${\rm Der}(N,M) \stackrel{\iota_1}{\hra} C^1(N,M)$, and we have 
 \BE\label{Autclassexpl} \omega^0_T(\hat{\partial}^1_N)(z) = [\iota_{1*}^{-1}{\partial}^1_{T}(\gamma)] \:.\EE


\N We are now ready to state our main result.\V

\begin{satz}\label{theorem} Let $G$ be the semidirect product of a normal subgroup $N$ and a subgroup $T$, and let $M$ be a $G$-module. Then sequence \REF{H2sequ} in the introduction is exact, where the maps are defined as follows.
For $d\in {\rm Der}(T,\Der)$ the class $\tau[d]$  is represented by the 2-cocycle $\beta_d\mapco G\times G \to M$, $\beta_d(nt,n't') = nd(t)({}^tn')$ for $n,n'\in N$, $t,t'\in T$, and the map $\phi$ is given by  $\phi=(\omega^0_T(\hat{\partial}^1_N)\,, \partial^0_{T*})$. 
\end{satz} 

The proof will occupy the rest of the paper.

\section{Automorphisms of group extensions}

We first recall some basic facts about group extensions and homomorphisms between them.

An extension of $G$ by $M$ is a short exact sequence of groups 
\[ {\cal E} \mapco \hspace{2mm} 0 \to M \mr{i} E \mr{\pi} G \to 1 \]
(which we will mostly write in the shorter form $M \Inj{i} E \Sur{\pi} G$) such that the given action of $G$ on $M$ coincides with the one induced by conjugation in $E$, i.e., ${}^ei(m) = i(\pi(e)m)$ for $e\in E$, $m\in M$. Two group extensions $ {\cal E}, {\cal E}'$ are said to be {\em congruent}\/ if there is a map (and hence an isomorphism) from $E$ to $E'$ commuting with the injections of $M$ and the projections to $G$. Congruence classes of extensions of $G$ by $M$ form an abelian group which is canonically isomorphic with $H^2(G,M)$, see \cite[IV.3]{ML}.
Finally, if $f\mapco \Gamma \to G$ is a homomorphism we denote by $f^*M$ the $\Gamma$-module which is $M$ as an abelian group endowed with the action of $\Gamma$ given by pulling back the action of $G$ via $f$.\V

\begin{prop}\label{Homext} Let ${\cal E}_k\mapco M \Inj{i_k} E_k \Sur{\pi_k} G$, $k=1,2$, be two group extensions of $G$ by $M$, and let $f\in {\rm Hom} (G_1, G_2)$ and $\alpha\in {\rm Hom}_{G_1}(M_1,f^*M_2)$. Then the diagram of unbroken arrows
\BE\label{Homextdia}
\begin{CD}\begin{matrix}
\ruled {\cal E}_2\mapco &M_2 & \Inj{i_2} & E_2 & \Sur{\pi_2} & G_2 \cr
 & \mapup{\alpha} & & \brokenup{\epsilon} \mapup{f} \cr
{\cal E}_1\mapco   & M_1 & \Inj{i_1} & E_1 & \Sur{\pi_1} & G_1 
\end{matrix}\end{CD}\EE
admits a filler $\epsilon$ (i.e.\ a group homomorphism from $E_1$ to $E_2$ rendering the diagram commutative) if and only if $\alpha_*[{\cal E}_1] - f^*[{\cal E}_2] = 0$ in $H^2(G_1,f^*M_2)$. Moreover, the group ${\rm Der}(G_1,f^*M_2)$ acts simply and transitively on the set $X_{(f,\alpha)}$ of all such fillers, by $(d+\epsilon)(e_1)
= i_2(d\pi_1(e_1))\epsilon(e_1)$ for $d\in {\rm Der}(G_1,f^*M_2)$, $\epsilon \in X_{(f,\alpha)}$ and $e_1\in E_1$.\hfbox
\end{prop}

 Now let ${\cal E} \mapco  M \Inj{i} E \Sur{\pi} G$ be an   extension of $G$ by $M$. Consider the following subgroups of  $\Aut(E)$ and of $\Aut(G)\times \Aut(M) $, resp.
 \[ \AutM(E) = \{\epsilon\in \Aut(E) \,|\, \epsilon(iM)=iM\} \]\[ \Aut(G,M) = \{(f,\alpha) \in \Aut(G)\times \Aut(M) \,|\, \mbox{$\forall (g,m)\in G\times M$, $\alpha(gm) = f(g)\alpha(m)$}\} \]
 A homomorphism $\rho \mapco \AutM(E) \to \Aut(G,M)$ is defined by $\rho(\epsilon) = (\epsilon_G,\epsilon_M)$ where $\epsilon_G$ and $\epsilon_M$ are induced by $\epsilon$. 
Moreover, the group $\Aut(G,M)$ acts on $(C^*(G,M),\partial_G^*)$ by automorphisms of complexes where 
\BE\label{AutMaction} (f,\alpha)\beta = \alpha_*((f^{-1})^{\times  n})^*\beta \EE
 for $\beta\in C^n(G,M)$. We write $\alpha_*f^*[\beta]$ for the induced action on $H^n(G,M)$.\V

 \begin{kor}\label{Autext} There is an exact sequence of groups
 \BE\label{Autextsequ} 0 \to {\rm Der}(G,M) \mr{(-)+id} \AutM(E) \mr{\rho} \Aut(G,M) \mr{\cal O} H^2(G,M) \EE
 where $(-)+id$ is a homomorphism and ${\cal O} = \partial^0_{{\rm Aut}(G,M)}[{\cal E}]$ is the inner derivation associated with the element $
 [{\cal E}]$ of the $\Aut(G,M)$-module $H^2(G,M)$. More explicitly, ${\cal O} $ is
 given by ${\cal O}(f,\alpha) = \alpha_*f^*[{\cal E}] - [{\cal E}]$.  \hfbox
 \end{kor}
 
 We also need to determine the cohomology class of the group extension 
 \BE\label{AutextKerO}  {\Aut({\cal E})} \mapco 0 \to {\rm Der}(G,M) \mr{(-)+id} \AutM(E) \mr{\rho} \Ker{\cal O} \to 1 \EE
 obtained from sequence \REF{Autextsequ} . It is easy to check that the action of $\Ker{\cal O}$ on ${\rm Der}(G,M) $ induced by conjugation in $\AutM(E)$ coincides with the natural action given by \REF{AutMaction}, i.e., $(f,\alpha)d = \alpha d f$.\V

 \begin{prop}\label{Autclass} The class of the extension \REF{AutextKerO} in $H^2(\Ker{\cal O},{\rm Der}(G,M))$ is given by the element $
 \omega^0_{{\rm Ker}({\cal O})}(\hat{\partial}^1_N)[{\cal E}]$.
 \end{prop}
 
\N More explicitly, let $\beta \mapco G\times G \to M$ be a 2-cocycle representing the extension ${\cal E}$. Then for $(f,\alpha) \in \Ker{\cal O}$ there exists $\gamma(f,\alpha) \in C^1(G,M)$ such that 
 \BE\label{defgamma} (f,\alpha) \beta - \beta = {\partial}^1_G(\gamma(f,\alpha) ) \:.\EE
 We thus get a map $\gamma \in C^1( \Ker{\cal O},C^1(G,M))$. Its image ${\partial}^1_{{\rm Ker}({\cal O})}(\gamma) \in Z^2(  \Ker{\cal O},$ $C^1(G,M))$ actually takes values in ${\rm Der}(G,M) \subset C^1(G,M)$, and we have 
 \BE\label{Autclassexpl} [{\Aut({\cal E})}] = [\iota_{1*}^{-1}{\partial}^1_{{\rm Ker}({\cal O})}(\gamma)] \:.\EE
 
 \proof Evaluating the maps in equation \REF{defgamma} on the couple $(f(g),f(g'))$ for $(g,g')\in G^2$  we get the following relation.
 \BE\label{gammaexpl} \alpha\beta(g,g') - \beta(f(g),f(g')) = f(g)\gamma(f,\alpha)(f(g')) - \gamma(f,\alpha)(f(gg')) + \gamma(f,\alpha)(f(g))\EE
Next we use $\beta$ to replace ${\cal E}$ by the congruent extension
\[   {\cal E'} \mapco 0 \to M \mr{i'} E' \mr{\pi'} G \to 1\]
where $E' = M\times G$ endowed with the group law $(m,g)(m',g') = (m+gm'+\beta(g,g')\,,gg' )$, $i'(m)=(m,1)$ and $\pi'(m,g) = g$. 
It is clear that extension  $\Aut({\cal E})$ is congruent with $\Aut({\cal E'})$, so we may replace it by the latter.
We construct a  normalized  set-theoretic section $s\mapco \Ker{\cal O} \to \AutM(E')$ of $\rho$, as follows. For $(f,\alpha) \in \Ker{\cal O}$ define a map $s(f,\alpha) \mapco E'\to E'$, $s(f,\alpha) (m,g) = (\alpha(m)+\gamma(f,\alpha)f(g), f(g))$. We must check that $s(f,\alpha) $ is a homomorphism; the third (and crucial)  equality in the following calculation follows from \REF{gammaexpl}.
\begin{eqnarray*}
& & \hspace{-10mm} s(f,\alpha) \Big((m,g)(m',g')\Big) \\
&=&  s(f,\alpha) \Big(m+gm'+\beta(g,g')\,,gg' \Big) \\
  &=& \Big( \alpha(m) + \alpha(gm') + \alpha\beta(g,g')+\gamma(f,\alpha) f(gg')\,,f(gg') \Big) \\
  &=&  \Big( \alpha(m) + f(g)\alpha(m') + \gamma(f,\alpha) f(g) + f(g)\gamma(f,\alpha) f(g')+\beta(f(g),f(g)') \,,f(gg') \Big) \\
  &=& \Big( \alpha(m)  + \gamma(f,\alpha) f(g) \,,f(g)\Big) \Big( \alpha(m') + \gamma(f,\alpha) f(g') \,,f(g')\Big) \\
  &=&  \Big(s(f,\alpha) (m,g)\Big)  \Big(s(f,\alpha) (m',g')\Big) 
\end{eqnarray*}
Moreover, the diagram 
\[\begin{matrix}
  \ruled M & \mr{i'}  & E' &  \mr{\pi'}  & G \cr
   \mapdown{\alpha} & &  \mapdown{s(f,\alpha)} & &  \mapdown{f}\cr
  M & \mr{i'}  & E' &  \mr{\pi'}  & G
  \end{matrix}\]
 obviously commutes, whence $s(f,\alpha) \in \AutM(E')$ and $\rho(s(f,\alpha)) =  (f,\alpha)$. Thus the extension $\Aut({\cal E'})$ is represented by the 2-cocycle $\beta' \in Z^2(\Ker{\cal O},{\rm Der}(G,M))$ defined by
 \[ \beta'\Big( (f,\alpha)\,,(f',\alpha') \Big) = ((-) + id)^{-1} \Big( s(f,\alpha) \circ s(f',\alpha') \circ s(ff',\alpha\alpha')^{-1} \Big) \]
 But 
 \begin{eqnarray*}
s(f,\alpha) \circ s(f',\alpha') (m,g) &=&  s(f,\alpha) \Big( \alpha'(m)+\gamma(f',\alpha')f'(g)\,, f'(g) \Big) \\
  &=&  \Big( \alpha \alpha'(m) + \alpha \gamma(f',\alpha')f'(g) + \gamma(f,\alpha)ff'(g)\,, ff'(g) \Big)
\end{eqnarray*}
while 
\begin{eqnarray*}
& & \hspace{-10mm}  \Big( \beta'\Big( (f,\alpha)\,,(f',\alpha') \Big) +id \Big) \circ s(ff',\alpha\alpha')(m,g) \\
  &=&  \Big( \beta'\Big( (f,\alpha)\,,(f',\alpha') \Big) +id \Big) \Big( \alpha\alpha'(m) + \gamma(ff',\alpha\alpha')ff'(g)\,, ff'(g) \Big) \\
  &=&  \Big( \beta'\Big( (f,\alpha)\,,(f',\alpha') \Big)ff'(g) \,,1\Big) \Big( \alpha\alpha'(m) + \gamma(ff',\alpha\alpha')ff'(g)\,, ff'(g) \Big) \\
  &=&  \Big( \beta'\Big( (f,\alpha)\,,(f',\alpha') \Big)ff'(g) +  \alpha\alpha'(m) + \gamma(ff',\alpha\alpha')ff'(g) \,,ff'(g) \Big)
\end{eqnarray*}
Thus
\[ \Big( \beta'\Big( (f,\alpha)\,,(f',\alpha') \Big)ff'(g)  =   \alpha \gamma(f',\alpha')f^{-1}ff'(g) - \gamma(ff',\alpha\alpha')ff'(g) + \gamma(f,\alpha)ff'(g) 
 \]
whence
\begin{eqnarray*}
 \beta'\Big( (f,\alpha)\,,(f',\alpha') \Big) &=& (f,\alpha)\gamma(f',\alpha') -  \gamma(ff',\alpha\alpha')+ \gamma(f,\alpha) \\
  &=&  \partial^1_{{\rm Ker}({\cal O})}(\gamma)\Big( (f,\alpha)\,,(f',\alpha') \Big)
\end{eqnarray*}  
This shows that  the map $\partial^1_{{\rm Ker}({\cal O})}(\gamma) \in Z^2({\rm Ker}({\cal O}),C^1(G,M))$ actually takes values in Der$(G,M)$ $\stackrel{\iota_1}{\hra} C^1(G,M)$, so $\beta' = (\iota_{1*})^{-1}\partial^1_{{\rm Ker}({\cal O})}(\gamma)$. Now recall that 
$ \omega^0_{{\rm Ker}({\cal O})}(\hat{\partial}^1_N) =
 \omega_0\omega_1$ where
 \[ H^0(\Ker{\cal O},H^2(G,M)) \mr{\omega_0} H^1(\Ker{\cal O},\Imm{\partial^1_G}) \mr{\omega_1} H^2(\Ker{\cal O},{\rm Der}(G,M) )\]
are the connecting homomorpisms  
induced by the short exact sequences of $\Ker{\cal O}$-modules
\[ 0 \to \Imm{\partial^1_G} \stackrel{\iota_0}{\hra} Z^2(G,M) \mr{q_2} H^2(G,M) \to 0 \]
\[ \ruled 0 \to  {\rm Der}(G,M) \stackrel{\iota_1}{\hra} C^1(G,M) \mr{\tilde{\partial}^1_G} \Imm{\partial^1_G} \to 0 \]
where $q_2$ is the canonical projection and $\tilde{\partial}^1_G$ is given by ${\partial}^1_G$. 
So $[\Aut({\cal E}')] = [\beta'] = [(\iota_{1*})^{-1}\partial^1_{{\rm Ker}({\cal O})}(\gamma)]=[
(\iota_{1*})^{-1}\partial^1_{{\rm Ker}({\cal O})}(\T{\partial}^1_{G*})^{-1}\T{\partial}^1_{G*}(\gamma)]=\omega_1[\T{\partial}^1_{G*}(\gamma)]$. But
\begin{eqnarray*}
[\T{\partial}^1_{G*}(\gamma)] &=& [(\iota_{0*})^{-1} {\partial}^1_{G*}(\gamma)] \quad \mbox{since $\partial^1_G = \iota_0\T{\partial}^1_{G*}$} \\
  &=& [(\iota_{0*})^{-1} \partial^0_{{\rm Ker}({\cal O})}(\beta)] \quad \mbox{by \REF{defgamma}} \\
  &=& [(\iota_{0*})^{-1} \partial^0_{{\rm Ker}({\cal O})}q_2^{-1} [{\cal E}]] \\
  &=& \omega_0 [{\cal E}]
\end{eqnarray*}
So $[\Aut({\cal E})] = [\Aut({\cal E}')] = \omega_1 \omega_0 [{\cal E}]$, as asserted.\hfbox

\section{Extensions of semidirect products}

From now on we suppose that $G=N \rtimes T$, writing $N\stackrel{\iota_N}{\hra} G$, $T\stackrel{\iota_T}{\hra} G$, and $\varphi\mapco T\to \Aut(N)$ for the action given by conjugation in $G$.\V


  \begin{defi}\label{real} Let $\underline{E}_N\mapco M \Inj{i_N} E_N \Sur{\pi_N} N$ and \ET\ be group extensions and 
$\T{\varphi} \mapco E_T \to  \AutM(E_N)$ be a homomorphism. We say that the triple $(\underline{E}_N,\underline{E}_T, \tphi)$ is {\em realizable}\/ if there exists an extension ${\cal E} \mapco M \Inj{i} E \Sur{\pi} G$ and a commutative diagram
\BE\label{ENEET}\begin{matrix}
\ruled  {\cal E}_N \mapco & M & \Inj{i_N} & E_N &\Sur{\pi_N} &N \cr
 &   \| & & \mapdown{i_1} & & \injdown{\iota_N} \cr
\ruled\ruleu   {\cal E}  \mapco & M & \mr{i} & E & \Sur{\pi} & G \cr
 &  \| & & \mapup{i_2} & & \injup{\iota_T}\cr
 {\cal E}_T \mapco &  M & \Inj{i_T} & E_T & \Sur{\pi_T} & T
\end{matrix} \EE
such that for $e_T\in E_T$, $e_N\in E_N$
\BE\label{Tphicond} \tphi(e_T)(e_N) = i_1^{-1}\Big( {}^{i_2(e_T)}i_1(e_N) \Big)\:. \EE
\end{defi}

\begin{prop}\label{realcond} A triple $(\underline{E}_N,\underline{E}_T, \tphi)$ is realizable  if and only if the following diagram  commutes.
\BE\label{Tphidia}\begin{matrix}
  \ruled\ruleu M & \Inj{i_T} & E_T & \Sur{\pi_T} & T \cr
  \mapdown{{} -\partial^0_N}  &  &  \mapdown{\tilde{\varphi}}  &  &  \mapdown{(\varphi,\psi_T)}  \cr
\ruled\ruleu   \Der   &   \mr{(-)+id}  &  \AutM(E_N)  &  \mr{\rho}  &  \Aut(N,M)
  \end{matrix}\EE
\end{prop}\V

The proof requires the following, certainly wellknown construction.\V

Let $\Gamma$ be a group. Recall that a $\Gamma$-group is a group $\Lambda$ endowed with a homomorphism $\alpha\mapco \Gamma \to \Aut(\Lambda)$; we write $\gamma\cdot \lambda = \alpha(\gamma)(\lambda)$. A homomorphism of $\Gamma$-groups is a homomorphism between $\Gamma$-groups which is $\Gamma$-equivariant.\V

\begin{prop}\label{amalgam} Let $K \ml{f} \Lambda \mr{g} \Gamma$ where $f$ is a homomorphism of $\Gamma$-groups and where $g$ is a precrossed module, i.e.\ a a homomorphism of $\Gamma$-groups where $\Gamma$ acts on itself by conjugation. Furthermore, suppose that ${g(\lambda)}\cdot\kappa ={}^{f(\lambda)}\kappa$ for all $(\lambda,\kappa)\in \Lambda \times K$. Then the amalgamated semi-direct product group $K\rtimes_{\Lambda} \Gamma=
K \rtimes  \Gamma/\{(f(\lambda),g(\lambda)^{-1}) \,|\,\lambda \in \Lambda\}$ is defined and has the following universal property: for any commutative diagram of homomorphisms of $\Gamma$-groups
\[\begin{matrix}
\Lambda & \mr{g} & \Gamma \cr
\mapdown{f}  &  &  \mapdown{h_{\Gamma}} \cr
K & \mr{h_K} & \Omega
\end{matrix}\]
where $\Gamma$ acts on $\Omega$ via $\gamma\cdot \omega = {}^{h_{\Gamma}(\gamma)}\omega$, there is a unique homomorphism $(h_K,h_{\Gamma}) \mapco $ $ K\rtimes_{\Lambda} \Gamma \to \Omega$ such that $(h_K,h_{\Gamma})qj_1  = h_K$ and $(h_K,h_{\Gamma})qj_2  = h_{\Gamma}$ where $j_1(\kappa) = (\kappa,1)$, $j_2(\gamma) = (1,\gamma)$, and $q\mapco K \times \Gamma \auf K\rtimes_{\Lambda} \Gamma$ is the canonical projection. Moreover, if $f$ resp.\ $g$ is injective, then so is $qj_2$ resp.\ $qj_1$.\hfbox
\end{prop}\V

\proofof{Proposition \ref{realcond}} Suppose that $(\underline{E}_N,\underline{E}_T, \tphi)$ is realizable. Then the right hand square of diagram 
\REF{Tphidia} commutes since the vertical maps both are induced by conjugation in $E$. The left hand square also commutes since
\begin{eqnarray*}
(-\partial^0_N(m) + id)(e_N)  &=&  i_N(-\partial^0_N(m)(\pi_Ne_N))e_N  \\
  &=& i_N(m-\pi_N(e_N)m)e_N \\
    &=& i_N(m) ({}^{e_N}i_N(m)^{-1}) e_N \\
      &=& i_N(m) e_N i_N(m)^{-1} \\
      &=& i_1^{-1}\Big({}^{ i_2(i_Tm)} i_1(e_N) \Big) \\
      &=& \tphi(i_Tm)(e_N)\:.
\end{eqnarray*}
Conversely, let $(\underline{E}_N,\underline{E}_T, \tphi)$  such that diagram \REF{Tphidia} commutes. Then the maps $E_N \ml{i_N} M \mr{i_T} E_T$ satisfy the hypothesis of Proposition \ref{amalgam} where $E_T$ acts on $E_N$ via $\tphi$. Indeed, $i_N$ is $E_T$-equivariant since by commutativity of the right hand square of diagram \REF{Tphidia},
\begin{eqnarray}
i_N(\pi_T(e_T)m)  &=&  i_N\psi_T(\pi e_T)(m) \nonumber \\
  &=& i_N \tphi(e_T)_M(m)  \nonumber \\
  &=&i_N\Big(i_N^{-1}\tphi(e_T)(i_Nm) \Big)  \nonumber \\
\label{iNETequ}  &=&  \tphi(e_T)(i_Nm)) \:,
\end{eqnarray}
while commutativity of left hand square of \REF{Tphidia} implies  
\begin{eqnarray*}
\tphi(i_Tm)(e_N)  &=&  i_N\Big( ({}-\partial^0_N m)(\pi_N e_N)\Big)e_N \\
  &=&  i_N(m-\pi_N(e_N)m)e_N \\
  &=& i_N(m)({}^{e_N}i_N(m)^{-1} )e_N \\
  &=& {}^{i_N(m)}e_N\:.
\end{eqnarray*}
Thus the amalgamated semidirect product $E = E_N \rtimes_{\tphi} E_T / \{(i_N(m),i_T(m)^{-1}) \,|\,m\in M\}$ is defined, as well as the homomorphism $\pi = (\iota_N\pi_N,\iota_T\pi_T)\mapco E \to G$; in fact, $\iota_N\pi_N$ is $E_T$-equivariant as 
\begin{eqnarray*}
\pi_N\tphi(e_T)(e_N) &=&  \tphi(e_T)_N(\pi_Ne_N)  \\
  &=&  \varphi(\pi_Te_T)(\pi_Ne_N)  \\
  &=&  {}^{\iota_T\pi_T(e_T)} \iota_N\pi_N(e_N) \:.
\end{eqnarray*}
Putting $i_k=qj_k$, $k=1,2$, and $i=i_1i_N=i_2i_T \mapco M \to E$ we obtain a commutative diagram \REF{ENEET}. The sequence 
\BE\label{} {\cal E}(\underline{E}_N,\underline{E}_T, \tphi) \mapco M \Inj{i} E \Sur{\pi} G \EE
is an extension: $i$ is injective as $i_1$ (see Proposition \ref{amalgam}) and $i_N$ are; $\pi$ is surjective as $N$ and $T$ are in its image; and if $\overline{(e_N,e_T)} \in \Ker{\pi}$, $\pi_Ne_N = \pi_Te_T = 1 $ as $N\cap T=\{1\}$, whence $(e_N,e_T) = (i_Nm_1,i_Tm_2)$ for some $m_1,m_2\in M$. But then\begin{eqnarray*}
q(e_N,e_T)  &\equiv&  q\Big( (i_Nm_1,i_Tm_2)(i_Nm_2,i_T(m_2)^{-1}) \Big)\\
  &=& q  (i_Nm_1+\tphi(i_Tm_2)(i_Nm_2)\,,1)  \\
  &=&  q(i_Nm_1 + i_Nm_2\,,1) \quad \mbox{by \REF{iNETequ}} \\
  &=& i(m_1+m_2)\:.
\end{eqnarray*}
It remains to check that  ${}^{q(e_N,e_T)}i(m) = i(\pi q(e_N,e_T)m)$. Indeed, in $E_N \rtimes_{\tphi} E_T$,
\begin{eqnarray*}
(e_N,e_T)j_1i_N(m) (e_N,e_T)^{-1}  &=& (e_N,e_T)(i_Nm,1)(e_N,e_T)^{-1} \\
   &=&  \Big( e_N \tphi(e_T)(i_Nm) \tphi({}^{e_T}1)(e_N^{-1}) \,, 1 \Big) \\
   &=&  \Big( e_N i_N(\varphi(\pi_T e_T)(m))e_N^{-1}\,,1 \Big) \quad \mbox{(right hand square of \REF{Tphidia})} \\
   &=& \Big( {}^{e_N} i_N(\pi_T (e_T) m)\,,1\Big) \\
   &=& \Big(i_N\Big(\pi_N(e_N)\pi_T(e_T)m\Big)\,,1\Big) \\
   &=& j_1i_N(\pi q(e_N,e_T)m)\;.
\end{eqnarray*}
Finally, condition \REF{Tphicond} is satisfied by definition of the semidirect product, whence $(\underline{E}_N,\underline{E}_T, \tphi)$ is realizable.\hfbox
\V

\begin{prop}\label{rm} (a) If $(\underline{E}_N,\underline{E}_T, \tphi)$ is realizable then the restricted extensions ${\cal E}(\underline{E}_N,\underline{E}_T, \tphi)_N$ and ${\cal E}(\underline{E}_N,\underline{E}_T, \tphi)_T$ are congruent with $\underline{E}_N$ and $\underline{E}_T$, resp.\V

(b) Any extension ${\cal E}$ of $G$ by $M$ is congruent with ${\cal E}({\cal E}_N,{\cal E}_T, \tphi)$ where $\tphi$ is given by \REF{Tphicond}.
\end{prop}\V

\proof Assertion (a) is immediate from diagram \REF{ENEET}. Now if ${\cal E} \mapco M \Inj{i} E$ $ \Sur{\pi} G$ is any extension then the inclusions of $E_N$ and $E_T$ into $E$ induce a surjective homomorphism $\xi \mapco E_N \rtimes_{\tphi} E_T \to E$ whose kernel is $\{(i_N(m)\,,i_T(m)^{-1}) \,|\, m\in M\}$, so $\xi$ induces the desired congruence from ${\cal E}({\cal E}_N,{\cal E}_T, \tphi)$ to ${\cal E}$.\hfbox\V

Now let $\underline{E}_N \mapco M \Inj{i_N} M\rtimes N \Sur{\pi_N} N$ and $\underline{E}_T \mapco M \Inj{i_T} M\rtimes T \Sur{\pi_T} T$ be the canonical split extensions. Then the bottom sequence in \REF{Tphidia} (which is exact by Proposition \ref{Autext}) is split by means of the canonical section $s\mapco \Aut(N,M) \to \AutM(M\times N)$, $s(f,\alpha) = \alpha \times f$. Hence we have a commutative diagram of homomorphisms with short exact rows
\BE\begin{matrix}
\ruled M  &  \Inj{i_T} &  M \rtimes T  &  \mr{\pi_T}  &  T \cr
\mapdown{{}-\partial^0_N}  &  &  \mapdown{{}-\partial^0_N \times (\varphi,\psi)^t}  &  & \mapdown{(\varphi,\psi)^t} \cr
\ruled\ruleu \Der  &  \Inj{}  &  \Der \times \Aut(N,M)  &  \Sur{}  & \Aut(N,M)  \cr
\|  &  &  \mapdown{\zeta}  &  &  \| \cr
\ruled \Der  & \Inj{(-)+id}  &  \AutM(M\rtimes N) & \Sur{\rho}  &  \Aut(N,M)
\end{matrix}\EE
where $\zeta(d,(f,\alpha)) = (d+id) \circ (\alpha \times f)$. Put $\tphi_0 = \zeta (\partial^0_N \times (\varphi,\psi)^t) \mapco M\rtimes T \to \AutM(M\rtimes N)$. Now let $d\in {\rm Der}(T,\Der)$ and $\tphi_d = d +\tphi_0$, see Proposition \ref{Homext}. Then
\begin{eqnarray}
\tphi_d(0,t)(m,n)  &=&  (d(t) + id) \circ \tphi_0(0,t)(m,n) \nonumber \\
  &=&  (d(t) + id) \circ (\psi(t) \times \varphi(t))(m,n) \nonumber\\
  &=&  (d(t) + id)(tm,{}^tn)  \nonumber\\
  &=&  i_Nd(t)\pi_N(tm,{}^tn)) (tm,{}^tn) \nonumber\\
  &=&  (d(t)({}^tn),1)(tm,{}^tn)\nonumber\\
  &=&  (d(t)({}^tn) +tm\,,{}^tn) \label{phid}
  \end{eqnarray}\vspace{-1mm}

\begin{prop}\label{Derprops} Let $d\in {\rm Der}(T,\Der)$. Then the 2-cochain $\beta_d\mapco G\times G \to M$, $\beta_d(nt,n't') = nd(t)({}^tn')$  is a 2-cocycle representing the extension ${\cal E}(\underline{E}_N,\underline{E}_T,\tphi_d)$. Moreover, the following properties are equivalent.\V

\N(1) $[\beta_d] = 0$ in $H^2(G,M)$;\V

\N(2) there exist derivations $D_N \mapco N \to M$ and $D_T \mapco T \to M$ such that $\beta_d$ is the coboundary of the function $D\mapco G\to M$ defined by $D(nt) = nD_T(t) + D_N(n)$;\V

\N(3) there exist derivations $D_N \mapco N \to M$ and $D_T \mapco T \to M$ such that $d =  \partial^0_T(D_N) - \partial^0_{N*}(D_T)$.\V
\end{prop}

\proof  Abbreviate $E_d = (M\rtimes N) \rtimes_{\tphi_d} (M\rtimes T) / \{(i_N(m),i_T(m)^{-1})(m)\,|\,m\in M\}$. Then a normalized set-theoretic section $\sigma$ of $\pi\mapco E \auf G$ is given by $\sigma(nt) = q((0,n),(0,t))$. Then
\begin{eqnarray*}
\sigma(nt) \sigma(n't')  &=&  q\Big( ((0,n),(0,t)) ((0,n'),(0,t')) \Big) \\
  &=&  q \Big( (0,n) \tphi_d(0,t)(0,n') \,, (0,t)(0,t') \Big) \\
  &=&  q \Big( (0,n)(d(t)({}^tn'),{}^tn') \,, (0,tt') \Big) \quad \mbox{by \REF{phid}}\\
  &=& q\Big( (nd(t)({}^tn')\,,n({}^tn'))\,,(0,tt') \Big) \\
  &=&  q\Big( (nd(t)({}^tn')\,,1)\,,(0,1) \Big) q\Big( (0,n({}^tn'))\,,(0,tt') \Big) \\
  &=&  i(nd(t)({}^tn')) \sigma(ntn't')
\end{eqnarray*}
Hence the 2-cocycle representing ${\cal E}(\underline{E}_N,\underline{E}_T,\tphi_d)$ associated to $\sigma$ is $\beta_d$. So it remains to prove the asserted equivalences. First note that the implication (2) $\Rightarrow$ (1) is plain, and that $\beta_d$ is the coboundary of a function $D\mapco G\to M$ iff $\forall (n,t),(n',t') \in N\times T$ one has
\BE\label{D}
D(ntn't') = ntD(n't') + D(nt) - nd(t)({}^tn') \:.\EE
Noting that for $t=1$ or $n'=1$ one has $d(t)({}^tn') =0$ we may take $t=t'=1$ or $n=n'=1$ in \REF{D} to see that the restriction of $D$ to $N$ and to $T$, denoted by $D_N$ and $D_T$, resp., are both derivations. Moreover, taking $t=n'=1$ in \REF{D} we get 
$D(nt') = nD(t') + D(n) = nD_T(t') + D_N(n) $, whence (1) implies (2). Now let $D_N\in \Der$ and $D_T\in {\rm Der}(T,M)$ and define $D\mapco G\to M$ by $D(nt) = nD_T(t)+ D_N(n) $. Then we have the following equivalences:
\begin{eqnarray*}
 & &  \mbox{$D$ satisfies \REF{D}}  \\
&\Longleftrightarrow& \rule{0mm}{12mm} \left\{ \begin{array}{ll} \ruled  & n({}^tn')D_T(tt') + D_N(n({}^tn'))  \cr = & nt(n'D_T(t') + D_N(n')) + nD_T(t) + D_N(n) - nd(t)({}^tn')\end{array}\right. \\
&\Longleftrightarrow& \rule{0mm}{10mm} \left\{ \begin{array}{ll}   \ruled  & ntn'D_T(t') + n({}^tn')D_T(t) +
nD_N({}^tn') + D_N(n)  \cr = & ntn'D_T(t') + ntD_N(n') + nD_T(t) + D_N(n) - nd(t)({}^tn')\end{array}\right. \\
&\Longleftrightarrow& \rule{0mm}{10mm} nd(t)({}^tn') = nt D_N(n') -nD_N({}^tn') + n(1-{}^tn' )D_T(t) \\
&\Longleftrightarrow& \rule{0mm}{10mm}  d(t)({}^tn') =  t D_N(n') - D_N({}^tn') - ({}^tn' -1)D_T(t)\:.
\end{eqnarray*}
Putting $n={}^tn'$ we see that $D$ satisfies \REF{D} iff $\forall (n,t)\in N\times T$, 
\begin{eqnarray*}
d(t)(n)  &=&  tD_N({}^{t^{-1}}n) - D_N(n) - (n-1)D_T(t) \\
  &=&  ( (t-1)D_N)(n) - \partial^0_N ( D_T(t))(n) \\
  &=&   \Big( \partial^0_T(D_N)(t) -   \partial^0_{N*}(D_T)(t) \Big)(n)\\
  &=&  \Big( \partial^0_T(D_N) -   \partial^0_{N*}(D_T)  \Big)(t)(n)\:.
\end{eqnarray*}
Hence (2) $\Leftrightarrow$ (3).\hfbox\V

To prove our main result it now suffices to assemble all the above propositions, as follows.\V

\proofof{theorem \ref{theorem}} Let $d\in {\rm Der}(T,\Der)$. If $d$ is inner, i.e.\ if $d =\partial^0_T(D_N)$ for some $D_N\in \Der$, we can take $D_T=0$ in Proposition \ref{Derprops} (3) to see that $[\beta_d] =0$, so $\tau$ is welldefined. Moreover,  $\tau[d]=0$ iff $[d] \in \Imm{\partial^0_{N*}}$, again by 
  Proposition \ref{Derprops}. Therefore sequence \REF{H2sequ} is exact in $H^1(T,\Der)$. To prove exactness in $H^2(G,M)$ first note that $res\circ \tau[d] = (res_N^G[{\cal E}(\underline{E}_N,\underline{E}_T,\tphi_d)]\,,res_T^G[{\cal E}(\underline{E}_N,\underline{E}_T,\tphi_d)]) = ([\underline{E}_N],[\underline{E}_T]) = (0,0)$ by 
Propositions \ref{Derprops}, \ref{rm}(a) and by construction of  ${\cal E}(\underline{E}_N,\underline{E}_T,\tphi_d)$. Now let $\cal E$ be some extension of $G$ by $M$ such that $res[{\cal E}]=0$. By Proposition \ref{rm}(b) $\cal E$ is congruent with ${\cal E}({\cal E}_N,{\cal E}_T, \tphi)$; as 
${\cal E}_N$ and ${\cal E}_T$ are split we may replace them by the canonical split extensions. The triple ${\cal E}({\cal E}_N,{\cal E}_T, \tphi)$ being realizable $\tphi$ fits into the commutative diagram \REF{Tphidia} by Proposition \ref{realcond}, so by Proposition \ref{Homext}, $\tphi = \tphi_d$ for some $d\in {\rm Der}(T,\Der)$. Thus $[{\cal E}] = [{\cal E}(\underline{E}_N,\underline{E}_T,\tphi_d)] = [\beta_d] = \tau[d]$ by Proposition \ref{Derprops}. Thus  sequence \REF{H2sequ} is exact in $H^2(G,M)$.
Finally, let $x=([\underline{E}_N],[\underline{E}_T])\in H^2(N,M)^T \times H^2(T,M)$. By Proposition \ref{realcond}, $x\in \Imm{res}$ iff there exists a homomorphism $\tphi \mapco E_T\to \AutM(E_N)$ fitting into the commutative diagram \REF{Tphidia}. Now $\Imm{(\varphi,\psi_T)} \subset \Ker{\cal O}$ since $[\underline{E}_N]$ is $T$-invariant, so by Proposition \ref{Homext} a filler $\tphi $ of \REF{Tphidia} exists iff $(\varphi,\psi_T)^*[\Aut(\underline{E}_N)] - ({}-\partial^0_N)_*[\underline{E}_T] = 0$ in $H^2(T,(\varphi,\psi_T)^*\Der) = H^2(T,\Der)$ by definition of the $T$-action on $C^*(N,M)$. But 
\begin{eqnarray*}
(\varphi,\psi_T)^*[\Aut(\underline{E}_N)] &=&  (\varphi,\psi_T)^*\omega_1\omega_0[\underline{E}_N)] \quad \mbox{by Proposition \ref{Autclass}}\\
  &=&   \omega_1\omega_0(\varphi,\psi_T)^*[\underline{E}_N)] \quad \mbox{by naturality of connecting maps}\\
  &=&   \omega_1\omega_0[\underline{E}_N)] 
\end{eqnarray*}
So $x\in \Ker{res}$ iff $\phi(x)=0$, which concludes the proof.\hfbox
\tw
\vspace{3mm}


\begin{thebibliography}{12}

\bibitem{Brown} K.\ S.\ Brown, Cohomology of groups, Springer GTM, Vol. 87, Springer-Verlag, New York,
Heidelberg,  Berlin, 1982.

\bibitem{Ga-Ma} H.\ Gaudier,   and R.\ Massy,   Une lin\'earisation de la notion de 2-cocycle, Pub. IRMA, Lille (1995), Vol.\ 37, N$^o$ II.


 \bibitem{HaRo} M.\ Hartl and C.\ Rousseau, Conceptual construction and proof of the fundamental exact sequence of group cohomology, in preparation.
 
 \bibitem{HueSS} J.\ Huebschmann, Automorphisms of group extensions and differentials in the Lyndon-Hochschild-Serre spectral sequence, J. of Algebra 72 (1981), 296-334.
 
 \bibitem{Hue8TS} J.\ Huebschmann, Group extensions, crossed pairs and an eight term exact sequence, J.\ Reine Angew.\ Math.\ 321 (1981), 150-172.

\bibitem{ML}  S.\ Mac Lane,   Homology, Springer Grundlehren, Vol.
114, Springer-Verlag Berlin-G\"ottingen-Heidelberg, 1963.

\bibitem{Ta}  K.-I.\ Tahara,
 On the second cohomology groups of semidirect products,    Math. Z. 129  (1972)   365-379.


\end{thebibliography}
